\documentclass[a4paper,12pt]{article}
\usepackage{amscd}
\usepackage{amsmath,amsfonts,amssymb,amscd}
\usepackage{indentfirst,graphicx,epsfig}
\usepackage{graphicx,psfrag}
\usepackage{ifpdf}
\usepackage{float}

\input{epsf}

\newtheorem{thm}{Theorem}

\newtheorem{lem}{Lemma}

\begin{document}
\title{\Large\bf On a question on graphs with\\ rainbow
connection number 2\footnote{Supported by NSFC No.11071130.}}
\author{\small Jiuying Dong, Xueliang Li\\
\small Center for Combinatorics and LPMC-TJKLC\\
\small Nankai University, Tianjin 300071, China\\
\small jiuyingdong@126.com, lxl@nankai.edu.cn}
\date{}
\maketitle

\begin{abstract}

For a connected graph $G$, the \emph{rainbow connection number
$rc(G)$} of a graph $G$ was introduced by Chartrand et al. In
``Chakraborty et al., Hardness and algorithms for rainbow
connection, J. Combin. Optim. 21(2011), 330--347", Chakraborty et
al. proved that for a graph $G$ with diameter 2, to determine
$rc(G)$ is NP-Complete, and they left 4 open questions at the end,
the last one of which is the following: Suppose that we are given a
graph $G$ for which we are told that $rc(G)=2$. Can we rainbow-color
it in polynomial time with $o(n)$ colors ? In this paper, we settle
down this question by showing a stronger result that for any graph
$G$ with $rc(G)=2$, we can rainbow-color $G$ in polynomial time by
at most 5 colors. \\[3mm]


\end{abstract}

\section{Introduction}

All graphs considered in this paper are simple, finite and
undirected. Undefined terminology and notation can be found in
\cite{Bondy-Murty}. Let $G$ be a graph, and $c: E(G) \rightarrow
\{1,2,\cdots,k\}, k\in N$ be an edge-coloring, where adjacent edges
may be colored the same. A graph $G$ is \emph{rainbow connected} if
for every pair of distinct vertices $u$ and $v$ of $G$, $G$ has a
$u-v$ path whose edges are colored with distinct colors. The minimum
number of colors required to make $G$ rainbow connected is called
the \emph{rainbow connection number} of $G$, denoted by $rc(G)$.
These concepts were introduced in \cite{Chartrand-Johns}. In
\cite{S-Ch} Chakraborty et al. proved that for a graph $G$ with
diameter 2, to determine $rc(G)$ is NP-Complete, and they left 4
open questions at the end, the last one of which is the following:

Suppose that we are given a graph $G$ for which we are told that
$rc(G)=2$. Can we rainbow-color it in polynomial time with $o(n)$
colors ? For the usual coloring problem, this version has been well
studied. It is known that if a graph is 3-colorable (in the usual
sense), then there is a polynomial time algorithm that colors it
with $\widetilde{O}(n^{3/14})$ colors \cite{BK}.

Li et al. \cite{LLL} and Dong and Li \cite{D-L} showed that if $G$
is a bridgeless graph with diameter 2, then $rc(G)\leq 5$, and Dong
and Li \cite{D-L} showed that the upper bound 5 is tight. At that
time we did not realize that we could solve the above open question.
Actually, from the proof of \cite{D-L} we can first deduce the
following result:

\begin{lem}\label{lem1}
For any bridgeless graph $G$ with diameter 2, we can rainbow-color
$G$ in polynomial time by at most 5 colors.
\end{lem}

Since a graph $G$ with a bridge that has $rc(G)=2$ must be composed
of a bridgeless graph $G'$ of radius 1 with a pendant edge attached
at the center vertex of $G'$, also from \cite{D-L} we can get the
following result:

\begin{lem}\label{lem2}
For any graph $G$ that is composed of a bridgeless graph $G'$ of
radius 1 with a pendant edge attached at the center vertex $G'$, we
can rainbow-color $G$ in polynomial time by at most 4 colors.
\end{lem}

Since a graph $G$ with $rc(G)=2$ is either a bridgeless graph with
diameter 2, or a graph with a bridge that is composed of a
bridgeless graph $G'$ of radius 1 with a pendant edge (the bridge)
attached at the center vertex of $G'$, as a consequence of the above
two lemmas, we settle down the last question in \cite{S-Ch} as
follows:

\begin{thm}\label{thm1}
For any graph $G$ with $rc(G)=2$, we can rainbow-color $G$ in
polynomial time by at most 5 colors.
\end{thm}

Before proceeding, we need some notation and terminology. For two
subsets $X$ and $Y$ of $V$, an $(X,Y)$-path is a path which connects
a vertex of $X$ and a vertex of $Y$, and whose internal vertices
belong to neither $X$ nor $Y$. We use $E[X,Y]$ to denote the set of
edges of $G$ with one end in $X$ and the other end in $Y$, and
$e(X,Y)=|E[X,Y]|$. The $k$-step open neighborhood of $X$ is defined
as $N^{k}(X)=\{v\in V(G)|d(v,X)=k, k\geq 0\}$. Let $N[S]=N(S)\cup
S$. For a connected graph $G$, the eccentricity of a vertex $v$ is
$ecc(v)=\max_{x\in V(G)}d_G(v,x)$. The radius of $G$ is
$rad(G)=\min_{x\in V(G)}ecc(x)$. The diameter of $G$ is $\max_{x\in
V(G)}ecc(x)$, denoted by $diam(G)$.

\section{Proof of the lemmas and theorem}

Throughout this section, the input graph $G$ is always bridgeless
with $n$ vertices, $m$ edges, and diameter 2. We say the colorings
of the following cycles containing $u$ to be $appropriate$: let
$C_3=uv_1v_2u$ be a 3-cycle where $v_1,v_2\in N^{1}(u)$, and let
$c(uv_1)=1, c(uv_2)=2, c(v_1v_2)=3(4)$; let $C_4=uv_1v_2v_3u$ be a
4-cycle where $v_1,v_3\in N^{1}(u), v_2\in N^{2}(u)$, and let
$c(uv_1)=1, c(uv_3)=2, c(v_1v_2)=3,c(v_3v_2)=4$; let
$C_5=uv_1v_2v_3v_4u$ be a 5-cycle where $v_1,v_4\in N^{1}(u),
v_2,v_3\in N^{2}(u)$, and let $c(uv_1)=1, c(uv_4)=2, c(v_1v_2)=3,
c(v_3v_4)=4, c(v_2v_3)=5$. We know that the shortest cycles passing
through $u$ are only the above mentioned $C_3, C_4$ or $C_5$.

Based on the proof of the main result of \cite{D-L}, we first give
an algorithm to rainbow-color a graph with diameter 2 by at most 5
colors.

{\bf Algorithm Rainbow-Color:}

\noindent {\bf Step 1}: Find the center vertex $u$ of $G$.\\
{\bf Step 2}:  Let $B_1,B_2,\cdots, B_b\subset N^{1}(u) $ satisfy
that for any $1\leq i\neq j \leq b$, $B_i\cap B_j=\emptyset$,
$B_i=N_{N^{1}(u)}[b_i]$, and $|B_i|\geq 2$. If $\bigcup_{i=1}^{b}
B_i=N^{1}(u)$, for $1\leq i\leq b$, let $c(ub_i)=1$, for any
$e \in E(u,B_i\setminus \{b_i\})$, let $c(e)=2$,
for any $e \in E(G[N^{1}(u)]))$, let
$c(e)=3$, and stop. Otherwise, proceed to
the next step.\\
{\bf Step 3}: Let $B_{b+1}\subset N^{1}(u)\setminus
\bigcup_{i=1}^{b} B_i$ be as large as possible that satisfies that
for any $w\in B_{b+1}, wb_i\not\in E(G)$, but $\exists
w'\in\bigcup_{i=1}^{b} (B_i\setminus\{b_i\})$ such that $ww'\in
E(G)$. If $\bigcup_{i=1}^{b+1} B_i=N^{1}(u)$, for $1\leq i\leq b$,
let $c(ub_i)=1$, for any $e\in E(u,B_i\setminus \{b_i\})$, let
$c(e)=2$, for any $e \in E(u, B_{b+1})$, let $c(e)=1$,
for any $e \in E(G[N^{1}(u)]))$, let
$c(e)=3$, and stop. Otherwise, let $B_{b+2}=N^{1}(u)\setminus \bigcup_{i=1}^{b+1} B_i$,
proceed to the next step.\\
{\bf Step 4}: Let $S=\bigcup_{i=1}^{b+1} B_i$.\\
While  $B_{b+2}\neq \emptyset$,\\
 For any any $v\in B_{b+2}$, we select a cycle $R$ such
that $R$ is a shortest cycle containing $uv$, and
 we further choose $R$ such that $R$ contains as
many vertices of $G\setminus S$ as possible.\\
If $V(R)\cap S=\{u\}$, then give $R$ an $appropriate$ coloring.
Otherwise,
give $R$ an $appropriate$ coloring according to the colors of colored edges of $R$.\\
Replace $S$ by $S\cup R$.\\
{\bf Step 5}: If $N^{2}(u)\subset S$, then for any $e\in
E(G[N^{1}(u)])$, let $c(e)=3$, and give the remaining
uncolored edges by a used color and stop. Otherwise, proceed to the next step.\\
{\bf Step 6}: Let $N^{1}(u)=X\cup Y$, where for any $x\in X$, $c(ux)=1$, for
any $y\in Y$, $c(uy)=2$.\\
{\bf Step 7}: Let all $S, T, Q \subseteq N^{2}(u)$ be as large as
possible which satisfy that
 for any $s\in S, E(s,X)\neq \emptyset$ but $E(s,Y)= \emptyset$;
for any $t\in T, E(t,Y)\neq \emptyset$ but $E(t,X)=
\emptyset$; for any $q\in Q, E(q,X)\neq \emptyset$ and $E(q,Y)\neq \emptyset$;
 for any $s\in S, E(s,T\cup Q)\neq\emptyset$; for any $t\in T,
E(t,S\cup Q)\neq\emptyset$. For any $e\in E(S\cup T\cup Q, X\cup
Y)$, give $e$ an $appropriate$ coloring.\\
{\bf Step 8}: If $N^{2}(u)= S\cup T\cup Q$, for any $e \in
E(G[N^{1}(u)])$, let $c(e)=3$; for any $e\in E(S, T\cup Q)$, let
$c(e)=5$, and we use a used color to color the remaining
uncolored edges and stop. Otherwise, proceed to the next stop.\\
{\bf Step 9}:  Let $ N^{2}(u)\setminus (S\cup T\cup Q)=P_1\cup P_2$
where for any $p_1\in P_1, e(p_1,X)=1$, for any
$p_2\in P_2, e(p_2,X)\geq 2$.\\
{\bf Step 10}: If  $P_1=\emptyset$, give a kind of coloring for the
remaining uncolored edges with 5 colors and stop.  Otherwise,
$P_1\neq\emptyset$.
If $|X|\geq 2$, go to Step 12. If $|X|=1$, go to the next step.\\
{\bf Step 11}: Let $D_1,D_2,\cdots, D_d\subset P$ satisfy that
 for $1\leq i\neq j \leq d, D_i\cap D_j=\emptyset,  D_i=N_P[d_i]$.\\
If $P=\bigcup_{i=1}^{d} d_i$, give a kind of coloring of the
remaining uncolored edges of $E(G)$ with 5 colors and stop.
Otherwise, let $D_{d+1}=P\setminus \bigcup_{i=1}^{d} d_i$, and give
a kind of coloring of the remaining uncolored edges of $E(G)$ with 5 colors and stop. \\
{\bf Step 12}: $\exists x\in X\setminus B_{b+2} $ such that
$E(x,P_1)=\emptyset$. Let $X_1\subseteq X$ be the set of all vertices which are adjacent
to the vertices of $P_1$. Let $X_2=X\setminus (X_1\cup B_{b+2})$.
Let $P_1'\subset P $ be the set of all vertices which are adjacent
to the vertices of $X_1$, and let $P_2'=P\setminus P_1'$ and give a kind of coloring of
the remaining uncolored edges of  $E(G)$ with 5 colors and stop.\\
{\bf Step 13}: For any $x\in X, E(x,P_1)\neq\emptyset$, give a kind
of coloring of the remaining uncolored edges of $E(G)$ with 5 colors
and stop.

Our algorithm is deduced from \cite{D-L}, and so we can
rainbow-color the input graph by at most 5 colors, which gives the
correctness of the algorithm. In the following we will examine the
time complexity of the algorithm.

At first, we know that each edge is colored only once by the
algorithm. Hence, the total effect of color assignments on the
algorithm's running time is $O(m)$. Thus, we can get that the
running time of Steps 5, 8, 10 and 13 is $O(m)$, respectively. The
time complexity for finding the center vertex in any connected graph
is $O(mn)$, because, it can be done by computing the eccentricity of
every vertex using a Breath First Search rooted at it. For any
vertex $v\in V(G)$, the running time for finding $N(v)$ is $O(n)$.
Hence the running time for Steps 2, 3, 7, 9, 11, and 12 is
$O(n^{2})$, respectively. For any vertex $v\in N^{1}(u)$, the time
complexity for finding $N_{N^{2}(u)}(v)$ is  $O(n)$; for any vertex
$v\in N^{2}(u)$, the time complexity for finding $N_{N^{2}(u)}(v)$
is $O(n)$; and for any vertex $v\in N^{2}(u)$, the time complexity
for finding $N_{N^{1}(u)}(v)$ is $O(n)$. Hence the total running
time for finding all the cycles $R$ and coloring all the cycles in
Step 4 is $O(n^{4})$. In Step 6, the time complexity is only
dependent on $N^{1}(u)$, hence the running time of Step 6 is $O(n)$.
Therefore, the total running time of the our algorithm is
$O(n^{4})$. The proof of Lemma 1 is now complete.

We know that any graph $G$ with a bridge that has $rc(G)=2$ must be
composed of a bridgeless graph $G'$ of radius 1 with a pendant edge
(the bridge) attached on the center vertex of $G'$. From the above
algorithm, we know that we can rainbow-color the bridgeless graph
$G'$ with radius 1 by at most 3 colors. For the pendant edge
attached on the center vertex $u$, we use another fresh color to
color the bridge. This gives $G$ a rainbow-coloring. Only Steps 1,2
and 3 of the algorithm are used, and the running time of them is
$O(n^{2})$, respectively. Hence, for any graph $G$ that is composed
of a bridgeless graph $G'$ of radius 1 with a pendant edge (the
bridge) attached on the center vertex of $G'$, we can rainbow-color
$G$ in polynomial time by at most 4 colors. The proof of Lemma 2 is
thus complete.

Finally, since a graph $G$ with $rc(G)=2$ is either a bridgeless
graph with diameter 2, or a graph with a bridge that is composed of
bridgeless graph $G'$ of radius 1 with a pendant edge (the bridge)
attached at the center vertex of $G'$, as a consequence of Lemmas 1
and 2, the proof of Theorem 1 is complete.

\end{document}